\documentclass{amsart}
\usepackage[british]{babel}
\usepackage{amsthm,amssymb,mathtools}
\usepackage{enumitem}
\usepackage[latin1]{inputenc}

\usepackage{url}
\usepackage{hyperref}

\newtheorem{theorem}{Theorem}

\theoremstyle{definition}
\newtheorem{definition}[theorem]{Definition}
\newtheorem{remark}[theorem]{Remark}
\newtheorem{example}[theorem]{Example}

\newcommand{\rca}{\mathbf{RCA}}
\newcommand{\aca}{\mathbf{ACA}}

\begin{document}


\title[Effective transfinite recursion in reverse mathematics]{What is effective transfinite recursion in reverse mathematics?}

\author{Anton Freund}
\address{Fachbereich Mathematik, Technische Universit\"at Darmstadt, Schlossgartenstr.~7, 64289~Darmstadt, Germany}
\email{freund@mathematik.tu-darmstadt.de}

\begin{abstract}
In the context of reverse mathematics, effective transfinite recursion refers to a principle that allows us to construct sequences of sets by recursion along arbitrary well orders, provided that each set is $\Delta^0_1$-definable relative to the previous stages of the recursion. It is known that this principle is provable in~$\mathbf{ACA}_0$. In the present note, we argue that a common formulation of effective transfinite recursion is too restrictive. We then propose a more liberal formulation, which appears very natural and is still provable in~$\mathbf{ACA}_0$.
\end{abstract}

\maketitle
{\let\thefootnote\relax\footnotetext{\copyright~2021 The Authors. \emph{Mathematical Logic Quarterly} published by Wiley-VCH GmbH.\\
This is the accepted version of a publication in \emph{Mathematical Logic Quarterly} 66:4 (2020)~479-483. Please cite the official journal publication, which is an open access article under the terms of the Creative Commons Attribution-NonCommercial-NoDerivs License.}

Effective transfinite recursion is a method from computability theory, which goes back to work of A.~Church, S.~Kleene and H.~Rogers (see~\cite[Section~I.3]{sacks-higher-recursion}). It does also have important applications in reverse mathematics (see e.\,g.~\cite{greenberg-montalban-ranked}). In~\cite{DFSW-effective-ramsey}, one can find an explicit formulation of effective transfinite recursion in reverse mathematics, together with a detailed proof in~$\aca_0$. It appears to be open whether the principle can be proved in $\rca_0$. 

In the present note, we propose a formulation of effective transfinite recursion that appears stronger than the one in~\cite{DFSW-effective-ramsey}. We will argue that this is ``the correct" formulation of the principle in reverse mathematics: It is very convenient for applications, seems to be as general as possible, and is still provable in~$\aca_0$.

To avoid misunderstanding, we point out that our aim is rather pragmatic: In recent work on fixed points of well ordering principles (see~\cite{freund-Pi12-induction}, in particular the proof of Theorem~5.11 and the paragraph before Definition~6.1), we have found it difficult to use effective transfinite recursion as formulated in~\cite{DFSW-effective-ramsey}. The present note is supposed to provide a formulation that is easier to apply. We do not know whether it is more general in a strict sense, i.\,e., whether our formulation and the one from~\cite{DFSW-effective-ramsey} can be separated over~$\rca_0$.

Let us recall the formulation of transfinite recursion in reverse mathematics. Working in second order arithmetic, we consider a well order~$X=(X,<_X)$. A family of sets $Y_x\subseteq\mathbb N$ indexed by $x\in X$ can be coded into the single set $Y=\{(x,n)\,|\,x\in X\text{ and }n\in Y_x\}$. More officially, any such family will be given as a set $Y\subseteq\mathbb N$, so that $n\in Y_x$ becomes an abbreviation for $(x,n)\in Y$. We write $Y^x=\{(x',n)\in Y\,|\,x'<_X x\}$ for the subfamily of sets with index below~$x\in X$. In a recursive construction, one defines $Y_x$ relative to $Y^x$, by stipulating $Y_x=\{n\in\mathbb N\,|\,\varphi(n,x,Y^x)\}$ for some formula $\varphi$, possibly with further parameters. To express that~$Y$ is the family defined by this recursive clause, we write
\begin{equation}\label{eq:recursive-hierarchy}
H_\varphi(X,Y)\quad:\Leftrightarrow\quad Y=\{(x,n)\,|\,x\in X\text{ and }\varphi(n,x,Y^x)\},
\end{equation}
adopting the notation from~\cite{DFSW-effective-ramsey}. From $H_\varphi(X,Y)$ and $H_\varphi(X,Y')$ we can infer $Y=Y'$, since a minimal element of the set $\{x\in X\,|\,Y_x\neq Y_x'\}$ would lead to a contradiction. Principles of transfinite recursion assert that there is a set~$Y$ with $H_\varphi(X,Y)$, for any well order~$X$ and each formula $\varphi$ from a certain class.  For effective transfinite recursion, the idea is that $\varphi$ should express a $\Delta^0_1$-property. Specifically, \cite[\mbox{Definition~6.5}]{DFSW-effective-ramsey} requires that $\varphi$ is a $\Sigma^0_1$-formula and that we have a further \mbox{$\Sigma^0_1$-formula}~$\psi$ with
\begin{equation}
\forall_{x\in X}\forall_{Z\subseteq\mathbb N}\forall_{n\in\mathbb N}(\varphi(n,x,Z)\leftrightarrow\neg\psi(n,x,Z)).
\end{equation}
However, the condition that this must hold for all~$Z\subseteq\mathbb N$ appears too strong: For many $\Delta^0_1$-definable notions, the $\Sigma^0_1$-definition and the $\Pi^0_1$-definition are only equivalent for objects from a certain class.

\begin{example}\label{ex:eff-rec-fast-growing}
For a total function $F:\mathbb N\to\mathbb N$, the iterates $F^{(m)}$ with~$m\in\mathbb N$ are given by $F^{(0)}(n):=n$ and $F^{(m+1)}(n):=F(F^{(m)}(n))$. The functions $F_k:\mathbb N\to\mathbb N$ at the finite stages of the fast growing hierarchy are defined by the clauses $F_0(n):=n+1$ and $F_{k+1}(n):=F_k^{(n)}(n)$. Let us discuss whether the hierarchy of functions $F_k$ can be constructed by effective recursion on~$k\in X=\mathbb N$, say over~$\aca_0$. In order to describe $F_{k+1}$ relative to~$F_k$, we need to express the definition of iterates: If $F:\mathbb N\to\mathbb N$ is total, then \mbox{$F^{(m)}(n)=n'$} is equivalent to the following: We have $n'=n_m$ for some\allowbreak\mbox{---or} equivalently every---sequence $\langle n_0,\dots,n_m\rangle$ with $n_0=n$ and $F(n_i)=n_{i+1}$ for all~$i<m$. Based on this observation, one readily constructs $\Sigma^0_1$-formulas $\varphi(n,x,Z)$ and $\psi(n,x,Z)$ with the following properties:
\begin{enumerate}
\item If we have $H_\varphi(\mathbb N,Y)$, then $Y_k$ is the graph of $F_k$.
\item The equivalence $\forall_{n\in\mathbb N}(\varphi(n,x,Z)\leftrightarrow\neg\psi(n,x,Z))$ holds for $x=0$, and for $x>0$ when $Z_{x-1}$ is the graph of a total function.
\end{enumerate}
However, the given $\Sigma^0_1$- and $\Pi^0_1$-definitions of $F^{(m)}(n)=n'$ need not agree when $F$ fails to be total. Hence the equivalence in~(2) may not hold for arbitrary $Z\subseteq\mathbb N$. In the present case, this is easy to repair: The functions $F_k$ are increasing, so that any sequence $\langle n_0,\dots,n_m\rangle$ that witnesses an equation $F_k^{(m)}(n)=n'$ in the aforementioned sense must satisfy~$n_0,\dots,n_m\leq n'$. This bound allows us to turn $\varphi$ into a $\Delta^0_0$-formula (in the original language or in a harmless extension, depending on our encoding of sequences). If we now set $\psi:=\neg\varphi$, then the equivalence in~(2) does clearly hold for all~$Z\subseteq\mathbb N$. Nevertheless, it would be preferable if we could apply effective transfinite recursion without proving monotonicity---in particular this would allows us to replace~$F_0$ by a base function that is not monotone.
\end{example}

As the example suggests, the equivalence $\forall_{n\in\mathbb N}(\varphi(n,x,Z)\leftrightarrow\neg\psi(n,x,Z))$ in the premise of effective transfinite recursion should not be required for arbitrary~$Z\subseteq\mathbb N$. In concrete applications, it will usually be clear which properties $Z$ needs to have. To find a general condition, we anticipate the construction of a set $Y$ with $H_\varphi(X,Y)$. Intuitively, the equivalence in question is only needed for $Z=Y^x$. To make this precise, we abbreviate
\begin{equation}
X\!\restriction\!x:=\{x'\in X\,|\,x'<_X x\}
\end{equation}
for each element $x$ of our well order~$X$. Note that $H_\varphi(X,Y)$ entails $H_\varphi(X\!\restriction\!x,Y^x)$, which determines $Y^x$ uniquely. This suggests the following:

\begin{definition}\label{def:eff-rec}
Effective transfinite recursion is the following principle, where $\varphi$ and $\psi$ range over $\Sigma^0_1$-formulas: If $X$ is a well order and we have
\begin{equation}
\forall_{x\in X}\forall_{Z\subseteq\mathbb N}(H_\varphi(X\!\restriction\!x,Z)\to\forall_{n\in\mathbb N}(\varphi(n,x,Z)\leftrightarrow\neg\psi(n,x,Z))),
\end{equation}
then there is a set $Y$ with $H_\varphi(X,Y)$.
\end{definition}

Let us now show how our version of effective transfinite recursion can be applied:

\begin{example}\label{ex:eff-rec-fast-growing2}
We take up the discussion from Example~\ref{ex:eff-rec-fast-growing}, still over~$\aca_0$. If we have $H_\varphi(\mathbb N\!\restriction\!x,Z)$, an induction over $y<x$ shows that $Z_y$ is the graph of a total function. As we have seen, this property ensures $\forall_{n\in\mathbb N}(\varphi(n,x,Z)\leftrightarrow\neg\psi(n,x,Z))$. Hence the hierarchy of functions $F_k$ with $k\in\mathbb N$ can be constructed by the principle of effective (transfinite) recursion, as specified in Definition~\ref{def:eff-rec}.
\end{example}

The given example of effective transfinite recursion is typical: We want to construct a hierarchy with certain properties---in this case, a hierarchy of total functions. The very same properties are supposed to ensure that the recursion step is given by a $\Delta^0_1$-definable relation. Since the hierarchy is constructed by recursion, the obvious way to establish the required properties is by induction. To apply effective transfinite recursion, we simply anticipate the inductive argument. On an informal level, this justifies our claim that Definition~\ref{def:eff-rec} provides the most general formulation of effective transfinite recursion in reverse mathematics.

Even though our approach sounds quite canonical, it does have some limitations: It is known that the relation $F_k(n)=n'$ is $\Delta^0_1$-definable over $\rca_0$ (even when we replace $k\in\mathbb N$ by ordinals $\alpha<\varepsilon_0$, see e.\,g.~\cite{sommer95}). Working in $\rca_0$, we can thus form the hierarchy of functions $F_k$ and prove basic facts about it. What we cannot show is that all functions~$F_k$ are total (since $k\mapsto F_k(k)$ is closely related to the Ackermann function). In view of this fact, the construction from Example~\ref{ex:eff-rec-fast-growing2} cannot be directly implemented in $\rca_0$ (even when all we want is a hierarchy of partial functions). To avoid this issue, we have used $\aca_0$ as base theory for our example. In fact, the extension of $\rca_0$ by $\Sigma^0_2$-induction would have done equally well. As mentioned before, it is open whether $\rca_0$ (or its extension by $\Sigma^0_2$-induction) proves the principle of effective transfinite recursion itself---but this is a different matter.

It is natural to ask whether our version of effective transfinite recursion entails the analogous principle with $H_{\neg\psi}$ at the place of $H_\varphi$ (and vice versa). We do not know the answer over~$\rca_0$, while $\aca_0$ allows us to formalize the following:

\begin{remark}
Consider $\Sigma^0_1$-formulas $\varphi$ and $\psi$. We assume that $X$ is a well order and that we have
\begin{equation}\label{eq:premise-ind-psi}
\forall_{x\in X}\forall_{Z\subseteq\mathbb N}(H_{\neg\psi}(X\!\restriction\!x,Z)\to\forall_{n\in\mathbb N}(\varphi(n,x,Z)\leftrightarrow\neg\psi(n,x,Z))).
\end{equation}
In order to apply effective transfinite recursion as formulated in Definition~\ref{def:eff-rec}, we show that $H_{\neg\psi}(X\!\restriction\!x,Z)$ follows from $H_\varphi(X\!\restriction\!x,Z)$. The latter entails \mbox{$H_\varphi(X\!\restriction\!y,Z^y)$} for all $y<_X x$. By transfinite induction over~$y<_X x$ we can establish the arithmetical statement $H_{\neg\psi}(X\!\restriction\!y,Z^y)$. Indeed, the induction hypothesis and~(\ref{eq:premise-ind-psi}) yield
\begin{equation}
\forall_{z<_X y}\forall_{n\in\mathbb N}(\varphi(n,z,Z^z)\leftrightarrow\neg\psi(n,z,Z^z)).
\end{equation}
In view of~(\ref{eq:recursive-hierarchy}), we can now infer $H_{\neg\psi}(X\!\restriction\!y,Z^y)$ from $H_\varphi(X\!\restriction\!y,Z^y)$. Having shown that $H_{\neg\psi}(X\!\restriction\!y,Z^y)$ holds for all $y<_X x$, one can use the same argument to deduce $H_{\neg\psi}(X\!\restriction\!x,Z)$ from $H_\varphi(X\!\restriction\!x,Z)$. As noted above, we can now apply the recursion principle from Definition~\ref{def:eff-rec}, which provides a set $Y$ with $H_\varphi(X,Y)$. Another transfinite induction yields $H_{\neg\psi}(X\!\restriction\! x,Y^x)$ for all $x\in X$, and then~$H_{\neg\psi}(X,Y)$.
\end{remark}

Let us now prove the following result, which is the main objective of our paper:

\begin{theorem}
Each instance of effective transfinite recursion (as formulated in Definition~\ref{def:eff-rec}) can be proved in~$\aca_0$.
\end{theorem}
As mentioned above, a detailed proof for an (apparently) weaker form of effective transfinite recursion has been given in~\cite[Proposition~6.6]{DFSW-effective-ramsey}. In the following, we recall this proof (including some notation from~\cite{DFSW-effective-ramsey}) and modify it in the relevant places. Let us mention that the recursion theorem can be used to give a different proof of the result (cf.~\cite[Section~I.3]{sacks-higher-recursion}).
\begin{proof}
For partial functions $f,g:\mathbb N\xrightarrow{p}\{0,1\}$ we write $f\preceq g$ to express that the graph of $f$ is contained in the graph of~$g$. To explain $f\preceq Z$ with $Z\subseteq\mathbb N$, we agree to identify $Z$ with its (total) characteristic function. From now on, the letters $f,g,h$ and $j$ are reserved for partial functions with finite domain, which we assume to be coded by natural numbers. Let us now fix $\Sigma^0_1$-formulas $\varphi$ and $\psi$, which determine an instance of effective transfinite recursion. The Kleene normal form theorem (see e.\,g.\ \cite[Theorem~II.2.7]{simpson09}) provides a $\Delta^0_0$-formula $\varphi_0$ such that $\aca_0$ (even $\rca_0$) proves
\begin{equation}\label{eq:varphi-finite}
\varphi(n,x,Z)\leftrightarrow\exists_{f\in\mathbb N}(f\preceq Z\land\varphi_0(n,x,f))\text{ and }\varphi_0(n,x,f)\land f\preceq g\to\varphi_0(n,x,g).
\end{equation}
Let $\psi_0$ be related to $\psi$ in the same way. Working in $\aca_0$, we assume that $X$ is a well order and that the premise
\begin{equation}
\forall_{x\in X}\forall_{Z\subseteq\mathbb N}(H_\varphi(X\!\restriction\!x,Z)\to\forall_{n\in\mathbb N}(\varphi(n,x,Z)\leftrightarrow\neg\psi(n,x,Z)))
\end{equation}
of effective transfinite recursion is satisfied. We want to construct finite approximations $f\preceq Y$ to a set $Y\subseteq\mathbb N$ with $H_\varphi(X,Y)$. In the process, we will also construct approximations to the sets~$Y^x$ for $x\in X$. To relate them, we agree to write $f\!\restriction\! x$ for the restriction of $f:\mathbb N\xrightarrow{p}\{0,1\}$ to arguments $(x',n)$ with $x'<_X x$. For $f\preceq g$, we say that $g$ is an $x$-extension of~$f$ if the following holds: When $g(x',n)$ is defined but $f(x',n)$ is not, then we have $x\leq_X x'$ or $x'\notin X$, as well as $g(x',n)=0$. As in the proof of~\cite[Proposition~6.6]{DFSW-effective-ramsey}, a finite partial function $f:\mathbb N\xrightarrow{p}\{0,1\}$ will be called an $x$-approximation if the following conditions are satisfied:
\begin{enumerate}[label=(\roman*)]
\item When $f(x',n)$ is defined with $x\leq_X x'$ or $x'\notin X$, then we have $f(x',n)=0$.
\item When we have $f(x',n)=1$ (resp.~$f(x',n)=0$) with $x'<_X x$, we have $\varphi_0(n,x',h)$ (resp.~$\psi_0(n,x',h)$) for some $x'$-extension $h$ of $f\!\restriction\! x'$.
\end{enumerate}
For an $x$-approximation $f$, one observes the following: If we have $x'<_X x$, then $f\!\restriction\!x'$ is an $x'$-approximation; and if $g$ is an $x$-extension of~$f$, then $g$ is an $x$-approximation as well. Consider the following two statements:
\begin{enumerate}[label=(\Roman*)]
\item If we have $\varphi_0(n,x,f)$ and $\psi_0(n,x,g)$ for some $n\in\mathbb N$, then $f$ and $g$ cannot both be $x$-approximations.
\item For any $x$-approximation $h$ and any $n\in\mathbb N$, there is an $x$-approximation $f$ with $h\preceq f$ that satisfies $\varphi_0(n,x,f)$ or $\psi_0(n,x,f)$.
\end{enumerate}
To avoid confusion, we note that~(II) will often be applied to the empty function~$h$, which is easily seen to be an $x$-approximation. In the proof of~\cite[Proposition~6.6]{DFSW-effective-ramsey}, statements (I) and (II) were established by (separate) inductions over~$x\in X$. Note that the required induction principle is available: Since the conjunction of (I) and~(II) is arithmetical, the set of counterexamples can be formed in $\aca_0$ (but not necessarily in~$\rca_0$). The proof from~\cite{DFSW-effective-ramsey} does not quite establish our version of effective transfinite recursion, since it relies on the assumtion that $\varphi(n,x,Z)\leftrightarrow\neg\psi(n,x,Z)$ holds for all values of the variables. However, we can adapt the argument to establish (I) and~(II) by simultaneous induction: Aiming at a contradiction, assume that~$x_0\in X$ is minimal such that (I) or (II) fails. We can then consider the $\Delta^0_1$-definable set
\begin{equation}
\begin{aligned}
Z&=\{(x,n)\,|\,x<_X x_0\text{ and }\varphi_0(n,x,f)\text{ for some $x$-approximation~$f$}\}=\\
{}&=\{(x,n)\,|\,x<_X x_0\text{ and }\psi_0(n,x,g)\text{ for no $x$-approximation~$g$}\}.
\end{aligned}
\end{equation}
In order to complete our induction, we must show that (I) and (II) hold for $x_0$. For this purpose, we will want to know that $\forall_{n\in\mathbb N}(\varphi(n,x_0,Z)\leftrightarrow\neg\psi(n,x_0,Z))$ holds for the set $Z$ that we have just defined. Due to the premise of effective transfinite recursion, the desired equivalence reduces to $H_\varphi(X\!\restriction\!x_0,Z)$. The latter will be established by another induction. As before, we write $Z^x=\{(x',n)\in Z\,|\, x'<_X x\}$ for $x<_X x_0$. Let us extend this notation by setting $Z^x:=Z$ for $x=x_0$. We now establish $H_\varphi(X\!\restriction\!x,Z^x)$ by induction over $x\leq_X x_0$. Due to the induction hypothesis, we can use the premise of effective transfinite recursion, which yields
\begin{equation}\label{eq:IH-Delta}
\forall_{y<_X x}\forall_{n\in\mathbb N}(\varphi(n,y,Z^y)\leftrightarrow\neg\psi(n,y,Z^y)).
\end{equation}
To establish $H_\varphi(X\!\restriction\!x,Z^x)$, we must prove that
\begin{equation}\label{eq:for-simult-ind}
(y,n)\in Z\quad\Leftrightarrow\quad\varphi(n,y,Z^y)
\end{equation}
holds for all $y<_X x$. As preparation, let us show that we have $f\preceq Z^y$ for any \mbox{$y$-}approximation~$f$. For later reference, we point out that the following goes through for any $y\leq_X x_0$. Given a value $f(y',n)=1$ (necessarily with $y'<_X y$), we obtain $\varphi_0(n,y',h)$ for some $y'$-extension $h$ of~$f\!\restriction\! y'$. As noted above, it follows that~$f\!\restriction\! y'$ and $h$ are $y'$-approximations. Hence $h$ witnesses $(y',n)\in Z$. Together with $y'<_X y$ we get $(y',n)\in Z^y$, as required. Now consider a value $f(y',n)=0$. If we have $y\leq_X y'$ or $y'\notin X$, then $(y',n)\notin Z^y$ is immediate. Hence we may assume $y'<_X y$. We then obtain $\psi_0(n,y',h)$ for some $y'$-approximation~$h$. This yields $(y',n)\notin Z$ and hence~$(y',n)\notin Z^y$. In order to deduce equivalence~(\ref{eq:for-simult-ind}), now for $y<_X x\leq_X x_0$, we first assume $(y,n)\in Z$. Then $\varphi_0(n,y,f)$ holds for some $y$-approximation~$f$. As we have seen, we get $f\preceq Z^y$. Due to equivalence~(\ref{eq:varphi-finite}), this yields $\varphi(n,y,Z^y)$. For the converse direction, we assume that $\varphi(n,y,Z^y)$ holds. By statement~(II) above, we may pick a $y$-approximation~$f$ with $\varphi_0(n,y,f)$ or $\psi_0(n,y,f)$. Aiming at a contradiction, we assume that $\psi_0(n,y,f)$ holds. Since $f$ is a $y$-approximation, we get $f\preceq Z^y$ and hence $\psi(n,y,Z^y)$. In view of~(\ref{eq:IH-Delta}), this contradicts the assumption that we have $\varphi(n,y,Z^y)$. So we must have $\varphi_0(n,y,f)$, which yields $(y,n)\in Z$. For~$x=x_0$, the result of the induction is $H_\varphi(X\!\restriction\!x_0,Z)$. Using the premise of effective transfinite recursion again, we get
\begin{equation}
\forall_{n\in\mathbb N}(\varphi(n,x_0,Z)\leftrightarrow\neg\psi(n,x_0,Z)),
\end{equation}
for~$Z$ as defined above. We now deduce that statements (I) and~(II) hold for~$x=x_0$. To establish~(I), assume that we have $\varphi_0(n,x_0,f)$ and $\psi_0(n,x_0,g)$. If $f$ and $g$ were $x_0$-approximations, then we would get $f,g\preceq Z$, as we have seen above. This would entail $\varphi(n,x_0,Z)$ and $\psi(n,x_0,Z)$, which contradicts the equivalence that we have just established. To establish~(II), we consider an $x_0$-approximation~$h$ and a number~$n\in\mathbb N$. Due to the last equivalence, we have $\varphi(n,x_0,Z)$ or $\psi(n,x_0,Z)$, which yields $\varphi_0(n,x_0,g)$ or $\psi_0(n,x_0,g)$ for some $g\preceq Z$. Note that $g$ and $h$ are compatible (i.\,e., have equal values whenever they are both defined) because of~\mbox{$h\preceq Z$}. By the monotonicity of $\varphi_0$ and $\psi_0$, it suffices to construct an $x_0$-approximation $f$ that satisfies $g\cup h\preceq f$. Let us enumerate the domain of $g\cup h$ as $\{(y_i,n_i)\,|\,i<k\}$. Following a hint by the referee, we simplify the argument by assuming~$y_0\leq_X\dots\leq_X y_{k-1}$. We will construct $x_0$-approximations $f_0\preceq\dots\preceq f_k$ such that the domain of $f_{i+1}$ contains $(y_i,n_i)$ but no elements $(x',m)$ with $y_i<_X x'$. Note that $g\cup h\preceq f_k$ is automatic due to~$f_k\preceq Z$. The construction proceeds exactly as in~\cite{DFSW-effective-ramsey}; we recall details in order to keep our presentation self-contained. To start the construction, let $f_0$ be the empty $x_0$-approximation. Assuming that $f_i$ is already defined, we construct $f_{i+1}$ by case distinction: We can put $f_{i+1}:=f_i$ when $(y_i,n_i)$ is already contained in the domain of~$f_i$; in the following we assume that this is not the case. If we have $x_0\leq_X y_i$ or $y_i\notin X$, it suffices to extend $f_i$ by the value $f_{i+1}(y_i,n_i):=0$. Now assume that we have $y_i<_X x_0$. By the induction hypothesis for~(II), we may pick a $y_i$-approximation $f$ with $f_i\!\restriction\!y_i\preceq f$ that satisfies $\varphi_0(n_i,y_i,f)$ or $\psi_0(n_i,y_i,f)$. Note that only one of these alternatives can hold, by the induction hypothesis for~(I). We can thus define $f_{i+1}$ as the extension of $(f\!\restriction\!y_i)\cup f_i$ by the value
\begin{equation}
f_{i+1}(y_i,n_i):=\begin{cases}
1 & \text{if $\varphi_0(n_i,y_i,f)$},\\
0 & \text{if $\psi_0(n_i,y_i,f)$}.
\end{cases}
\end{equation}
It is clear that $f_{i+1}$ satisfies condition~(i) from the definition of $x_0$-approximation. To verify condition~(ii), we consider a value $f_{i+1}(x',m)=1$ with $x'<_X x_0$ (the argument for $f_{i+1}(x',m)=0$ is completely parallel). Given that $(x',m)$ lies in the domain of~$f_{i+1}$, we even have $x'\leq_X y_i$ by construction. Let us consider $x'<_X y_i$ first. We then have $f(x',m)=1$. Since $f$ is a $y_i$-approximation, we get $\varphi_0(m,x',j)$ for some $x'$-extension~$j$ of~$f\!\restriction\! x'$. In view of $x'<_X y_i$ we have $f\!\restriction\! x'=f_{i+1}\!\restriction\!x'$, so that the same $j$ witnesses condition~(ii) for~$f_{i+1}$. Next, we consider $(x',m)=(y_i,n_i)$. Since we have assumed $f_{i+1}(x',m)=1$, we must have $\varphi_0(n_i,y_i,f)$. Let us observe that $f$ is a $y_i$-extension of $f_{i+1}\!\restriction\!y_i=f\!\restriction\!y_i$ (note that $f(x,k)=0$ holds if $y_i\leq_X x$ or $x\notin X$, since $f$ is a $y_i$-approximation). Hence $f$ itself witnesses condition~(ii) for the value~$f_{i+1}(y_i,n_i)$. Finally, assume we have $x'=y_i$ but $m\neq n_i$. In this case we get $f_i(x',m)=f_{i+1}(x',m)=1$. As $f_i$ is an $x_0$-approximation, we obtain $\varphi_0(m,x',j)$ for some $x'$-extension $j$ of~$f_i\!\restriction\!x'$. Recall that we have $f_i\!\restriction\!x'=f_i\!\restriction\!y_i\preceq f$, so that $f\!\restriction\!x'$ and $j$ are compatible. One can check that $j':=(f\!\restriction\!x')\cup j$ is an $x'$-extension of $f_{i+1}\!\restriction\!x'$. In view of $j\preceq j'$, we still have $\varphi_0(m,x',j')$. Hence~$j'$ witnesses condition~(ii) for the value~$f_{i+1}(x',m)$. This completes our inductive proof that statements~(I) and~(II) hold for all~$x\in X$. We can now set
\begin{equation}
\begin{aligned}
Y&=\{(x,n)\,|\,x\in X\text{ and }\varphi_0(n,x,f)\text{ for some $x$-approximation~$f$}\}=\\
{}&=\{(x,n)\,|\,x\in X\text{ and }\psi_0(n,x,g)\text{ for no $x$-approximation~$g$}\}.
\end{aligned}
\end{equation}
Consider the well order~$X\cup\{\top\}$ with a new maximum element~$\top$. Write \mbox{$X\!\restriction\!\top:=X$} and $Y^\top:=Y$. Just as above, an induction over $x\in X\cup\{\top\}$ yields~$H_\varphi(X\!\restriction\!x,Y^x)$. For $x=\top$ we obtain $H_\varphi(X,Y)$, as required by effective transfinite recursion.
\end{proof}

Our aim was to provide a convenient and canonical formulation of effective transfinite recursion in the setting of reverse mathematics. In our opinion, the discussion after Example~\ref{ex:eff-rec-fast-growing2} shows that we have achieved this aim. Since our objective was rather pragmatic, we have not considered the following questions, even though they are certainly interesting: Can~$\rca_0$ prove some formulation of effective transfinite recursion? Can it prove that the formulation from Definition~\ref{def:eff-rec} is equivalent to the one from~\cite{DFSW-effective-ramsey}?

\textbf{Acknowledgements.} I want to thank Linda Brown Westrick for her feedback on a first version of this paper.

\bibliographystyle{amsplain}
\bibliography{Effective-transf-rec_RM}

\end{document}